\documentclass[12pt]{amsart}
\usepackage{amssymb}
\usepackage{float}
\usepackage{graphicx}
\usepackage{textcomp}
 \usepackage{xcolor}
\usepackage[colorlinks=true, linkcolor=red, citecolor=blue]{hyperref}
\theoremstyle{plain}

\numberwithin{equation}{section}

\begin{document}
\title{Une $q$-d\'eformation de la  transformation  de Bargmann polyanalytique}
\author{S. Arjika$^{\dag,\ddag }$, O. El Moize$^{\flat }$ et Z. Mouayn$^{\ast }$}
\maketitle
\begin{center}
\begin{scriptsize}
$^{\dag }$  Department of Mathematics and Computer
Sciences, Faculty of Sciences and Technics,\vspace*{-0.5em}\\ University of Agadez, BP. 199, Agadez, Rep. of Niger\vspace*{0.2mm}\\
$^{\ddag }$ International Chair in Mathematical Physics and
Applications {\scriptsize(ICMPA-UNESCO Chair)},\\
 University of Abomey-Calavi, 072 BP.50,
Cotonou, Rep. of Benin\\
${}^{\flat} $ 46, Lot EL YOUSSR 2, Berrechid, Morocco\vspace*{0.2mm}\\
$^{\ast }$ Department of Mathematics, Faculty of Sciences
and Technics (M'Ghila),\vspace*{-0.5em}\\ BP. 523, B\'{e}ni Mellal, Morocco
\end{scriptsize}
\end{center}
\thanks{}
\begin{abstract}
\scriptsize{
We introduce a $q$-analog of the polyanalytic Bargmann transform on $\mathbb{C}$.}
\end{abstract}
{\scriptsize $\quad\:$	R\'ESUM\'E. Nous  introduisons  une version $q$-deform\'{e}e de la transformation de Bargmann
polyanalytique sur $\mathbb{C}$.}

\section{Introduction et \'enonc\'e des r\'esultats }
Dans \cite{B},  Bargmann  avait introduit une transformation bien c\'el\`ebre qui applique isom\'etriq-uement l'espace $L^2(\mathbb{R})$ sur l'espace de Fock des fonctions enti\`eres de carr\'e int\'egrables par rapport \`a la mesure Gaussienne $e^{-z\bar{z}}d\lambda$ o\`{u} $d\lambda$ d\'{e}signe la mesure de Lebesgue sur $\mathbb{C}$. \'{E}tant fortement li\'{e}e au groupe de Heisenberg cette transformation peut \^{e}tre vue comme une transformation de Fourier avec  fen\^{e}trage  \cite{BH}. D'o\`{u} le r\^{o}le important qu'elle joue en traitement du signal et dans l'analyse harmonique sur l'espace des phases \cite{Folland}.\medskip

 Il est aussi possible d'intrepr\'eter le noyau  de cette transformation en termes des \'etats coh\'erents \cite{AGA} associ\'es \`a l'Hamiltonien  de l'oscillateur harmonique dont les \'etats quantiques appartiennent \`a $L^2(\mathbb{R})$. Plus pr\'{e}cisement, un \'{e}tat coh\'{e}rent est repr\'{e}sent\'{e} par une fonction d'onde normalis\'{e}e que l'on d\'{e}finit \`{a} l'aide d'une superposition assez particuli\`{e}re de fonctions propres de l'Hamiltonien. C'est un fait bien connu que les fonctions propres de l'oscillateur harmonique sont donn\'{e}es par les fonctions d'Hermite. Il se trouve aussi que dans la superposition de ces \'{e}tats qui d\'{e}finit l'\'{e}tat coh\'{e}rent $|z\rangle$ index\'e par le nombre $z=q+ip$ appartenant \`a l'espace des phases $\mathbb{C}$, les coefficients apparaissent  sous la forme
\begin{equation}
h_j(z):=\frac{z^j}{\sqrt{j!}}, \ j=0,1,2,...\ .
\end{equation}
qui n'est autre que celle des \'el\'ements de la base de l'espace de Fock qui sera d\'esign\'e par $\mathfrak{F}(\mathbb{C})$. En notant $\mathcal{B}$ cette transformation, l'image d'une fonction arbitraire $f\in L^2(\mathbb{R})$ s'\'{e}crit 
\begin{equation}\label{barcl}
\mathcal{B}[f](z):=\pi^{-\tfrac{1}{4}}\int_{\mathbb{R}}e^{-\tfrac{1}{2}z^2-\tfrac{1}{2}\xi^2+\sqrt{2}\xi z}f(\xi)d\xi,\quad z\in\mathbb{C}.
\end{equation}
Par ailleurs, il a \'et\'e d\'emontr\'e  \cite{AIM} que l'espace  $\mathfrak{F}(\mathbb{C})$ co\"{i}ncide avec le noyau
\begin{equation}
\mathcal{A}_0(\mathbb{C}):=\{\varphi \in L^2(\mathbb{C},e^{-z\bar{z}}d\lambda),\tilde{\Delta}\varphi=0\}
\end{equation}
de l'op\'erateur diff\'erentiel du  second ordre 
\begin{equation}
\tilde{\Delta}:=-\frac{\partial^2}{\partial z\partial \bar{z}}+\bar{z}\frac{\partial}{\partial \bar{z}}.
\end{equation}
Ce dernier, suppos\'{e} agir sur l'espace de Hilbert $\mathtt{H}:=L^2\left(\mathbb{C},e^{-z\bar{z}}d\lambda\right)$, peut suite \`a un entrelacement unitaire (\textit{ground state transformation}), para\^{\i}tre sous la forme de l'op\'erateur de Schr\"{o}dinger associ\'e au mouvement projet\'e  sur le plan  $\mathbb{R}^2$ d'une particule \`a spin non nul, charg\'ee et plong\'ee dans un champ magn\'etique  uniforme et normal au plan. Le spectre de $\tilde{\Delta}$ relativement \`a $\mathtt{H}$ est constitu\'e de valeurs propres $\epsilon_m:=m\in \mathbb{Z}_+$, chacune \'{e}tant de multiplicit\'e infinie, appel\'ees niveaux de Landau Euclidiens. Notons qu'\`{a} chaque valeur propre est associ\'e un espace propre
\begin{equation}
\mathcal{A}_m(\mathbb{C}):=\{\varphi \in\mathtt{H
	},\tilde{\Delta}\varphi=\epsilon_m\varphi\}
\end{equation}
aussi appel\'e  \textquotedblleft\textit{true-polyanalytic Bargmann space}\textquotedblright $\:$ dans \cite{Vasilevski,AF}  dont une base orthonormale est donn\'ee par les fonctions
\begin{equation}
\label{truebase}
h_j^m(z):=(-1)^{m \wedge j}\left(m!j!\right)^{-1/2}  (m\wedge j)!|z|^{|m-j|}e^{-i(m-j)arg(z)}L_{m\wedge j}^{(|m-j|)}(z\bar{z}),\:z\in \mathbb{C}
\end{equation}
en termes des polyn\^omes de Laguerre $L^{(\alpha)}_n(.)$ (\cite{KS}, p.47)  o\`u $m\wedge j=min(m,j).$\medskip

Il fut donc naturel de faire jouer aux \'el\'ements $h_j^m(z)$ le m\^eme r\^ole que celui  des $h_j(z)$ \`a savoir  des coefficients dans une nouvelle superposition des fonctions propres de l'oscillateur harmonique. Les \'etats coh\'erents qui en r\'esultaient ont fourni une transformation de Bargmann g\'en\'eralis\'ee  not\'ee $\mathcal{B}_m:L^2(\mathbb{R})\rightarrow \mathcal{A}_m(\mathbb{C})$ ayant pour expression \cite{Mouayn1}:
\begin{equation}\label{tr1}
\mathcal{B}_m[f](z)=(-1)^m(2^mm!\sqrt{\pi})^{-\tfrac{1}{2}}\int_{\mathbb{R}}e^{-\tfrac{1}{2}z^2-\tfrac{1}{2}\xi^2+\sqrt{2}\xi z}H_m\left( \xi-\frac{z+\bar{z}}{2}\right)f(\xi)d\xi
\end{equation}
o\`{u} $H_m(.)$ est le polyn\^ome d'Hermite (\cite{KS}, p.59).  Plus d'informations sur $\mathcal{B}_m$ se trouvent dans \cite{AF} et les r\'ef\'erences qui  y figurent.\medskip

Notons que les coefficients \eqref{truebase} se laissent aussi s'exprimer \`a l'aide des polyn\^omes d'Hermite complexes \`a deux dimensions, not\'es $H_{r,s}(z_1,z_2)$, qui furent introduits par It\^o  \cite{IK} dans le contexte des processus de Markov complexes. Pr\'ecis\'ement, $ \left(m!j!\right)^{1/2}h_j^m(z)=H_{m,j}(z,\bar{z})$ o\`u
\begin{equation}\label{2dhermite}
H_{r,s}(z,w)=\displaystyle\sum_{k=0}^{r\wedge s} (-1)^k k!\binom{r}{k}\binom{s}{k} z^{r-k}w^{s-k},\quad r,s=0,1,2,...\ .
\end{equation}

Tout r\'ecemment, Ismail et Zhang \cite{IZ} ont introduit une version $q$-d\'eform\'ee des polyn\^omes $H_{r,s}(z,w)$ dont l'expression est
\begin{equation}
H_{r,s}(z,w|q)=\displaystyle\sum_{k=0}^{r\wedge s} \begin{bmatrix} r\\k \end{bmatrix}_{q} \begin{bmatrix} s\\k \end{bmatrix}_{q} (-1)^k q^{\binom{k}{2}}(q;q)_k z^{r-k}w^{s-k}, \ z,w\in\mathbb{C}
\end{equation}
o\`u
\begin{equation}
\begin{bmatrix} n\\ k \end{bmatrix}_q:=\frac{(q;q)_{n}}{%
	(q;q)_{n-k}(q;q)_{k}},\qquad (a;q)_{n}:=\medskip \prod\limits_{k=0}^{n-1}\left( 1-aq^{k}\right)
,\text{ } \: k=0,1,\cdots,n.
\end{equation}
Ces polyn\^{o}mes s'\'{e}crivent aussi 
\begin{equation}
\label{1.11}
H_{r,s}(z,w|q)=(-1)^{r\wedge s} \frac{(q;q)_{r \vee s}}{(q;q)_{|r-s|}}q^{\binom{r\wedge s}{2}}|z|^{|r-s|}e^{-i(r-s)arg( z)}P_{r\wedge s}\left( zw;q^{|r-s|}|q\right)
\end{equation}
\`a l'aide des polyn\^omes de Wall $P_n(\cdot,a|q)$ (\cite{KS}, p.107) o\`{u} $r\vee s=\max(r,s).$\medskip

 Ce nouveau mat\'eriel nous conduit donc \`a proposer une version $q$-d\'{e}form\'{e}e de la transformation \eqref{tr1}. Le noyau d'une telle transformation sera obtenu en superposant des $q$-d\'{e}form\'{e}es des fonctions d'Hermite \`{a} l'aide des $q$-d\'{e}form\'{e}s des coefficients $h_j^{m}(z)$, que nous proposerons de la forme 
\begin{equation}\label{Hmj}
h_j^{m,q}(z):=\frac{(-1)^{m\wedge j} (q;q)_{m \vee j}q^{\binom{m\wedge j}{2}}\sqrt{1-q}^{|m-j|}|z|^{|m-j|}e^{-i(m-j)\theta}}{(q;q)_{|m-j|}\sqrt{q^{mj}(q;q)_m(q;q)_j}}P_{m\wedge j}\left( (1-q)z\bar{z};q^{|m-j|}|q\right)
\end{equation}
que l'on se procure \`{a} partir de $(\ref{1.11})$ en rempla\c{c}ant $w$ par $\bar{z}$ et en normalisant le polyn\^{o}me obtenu. Et pour s'assurer de la consistance d'un tel choix de coefficients, il suffira de revenir \`a la d\'efinition du polyn\^ome de Wall
\begin{equation}\label{Wall1}
P_n(x;a|q) := \setlength\arraycolsep{1pt}
{}_2 \phi_1\left(\begin{matrix}q^{-n},0 \\ aq \end{matrix}\left|q;qx\right.\right)
\end{equation}
\`a l'aide de la fonction $q$-hyperg\'eom\'etrique ${}_2 \phi_1$ dont le comportement
\begin{equation}
\displaystyle\lim_{q\rightarrow 1} \setlength\arraycolsep{1pt}
{}_2 \phi_1\left(\begin{matrix}q^{-n},0 \\ q^{\alpha+1} \end{matrix}\left|q;q(1-q)x\right.\right)= \setlength\arraycolsep{1pt}
{}_1 F_1\left(\begin{matrix}-n \\ \alpha+1 \end{matrix}\left|x\right.\right)=\frac{n!}{(\alpha+1)_n}L_n^{(\alpha)}(x)
\end{equation}
permettra de r\'ecup\'erer le polyn\^ome de Laguerre  en tenant compte que ce dernier s'exprime aussi en terme de la  fonction hyperg\'eom\'etrique ${}_1F_1$. Effectivement, un calcul d\'etaill\'e permet d'aboutir \`a
$
\lim_{q\rightarrow 1} h_j^{m,q}(z)=h_j^m(z).
$

A pr\'{e}sent, nous choisissons comme $q$-analogues des \'etats propres de l'oscillateur harmonique les fonctions d\'efinies par 
\begin{equation}\label{fjq123}
\varphi_j^q(\xi):=\sqrt{\frac{\sqrt{2}\,\omega_q(\sqrt{2}\xi)}{(q;q)_j}}H_j\left(\sqrt{\frac{1-q}{2}}\xi|q\right)
\end{equation}
en termes des polyn\^omes $q$-Hermite continus $H_n(\cdot|q)$ (\cite{MDG}, p.381),  o\`u 
\begin{equation}
\omega_q(u):=\frac{(q;q)_{\infty}\sqrt{1-q}}{4\pi\sqrt{1-(1-q)u^2/4}}\displaystyle\prod_{k\geq 0}(1+(2-u^2(1-q))q^k+q^{2k}).
\end{equation}
Les fonctions que nous avons consid\'er\'ees dans \eqref{fjq123} forment un syst\`eme orthonormal sur l'intervalle  $\mathcal{I}_q:=]\frac{-\sqrt{2}}{\sqrt{1-q}},\frac{\sqrt{2}}{\sqrt{1-q}}[$. C'est-\`a-dire,
\begin{equation}
\displaystyle\int_{\mathcal{I}_q} \varphi_j^q(\xi)\varphi_k^q(\xi)d\xi=\delta_{jk}.
\end{equation} 

Cela nous permettra de d\'{e}finir un nouveau \'etat coh\'erent via la superposition
\begin{equation}\label{CS1}
|z,m,q\rangle=(\mathcal{N}_{m,q}(z\bar{z}))^{-\tfrac{1}{2}}\sum_{j\geq 0}\overline{h_j^{m,q}(z)}|\varphi_j^q\rangle,
\end{equation}
o\`u $\left(\mathcal{N}_{m,q}(z\bar{z})\right)^{\tfrac{-1}{2}}$ est un facteur de normalisation que l'on peut d\'eduire par calcul \`a partir de la relation \eqref{kern1} ci dessous. Pr\'ecis\'ement
\begin{equation}\label{normfac}
\mathcal{N}_{m,q}(z\bar{z})=\frac{q^{-m}(q^{1-m}(1-q)z\bar{z};q)_m }{(q^{-m}(1-q)z\bar{z};q)_{\infty}}.
\end{equation}
Cependant, la condition de finitude du facteur \eqref{normfac} contraint la variable d'indexation $z$ \`a ne pas quitter le domaine 
\begin{equation}
\mathbb{C}_{q,m}:=\{z\in \mathbb{C},\,(1-q)z\bar{z}<q^m\} 
\end{equation}
chose que l'on peut d\'eduire en exigeant \`a la quantit\'e $(q^{-m}(1-q)z\bar{z};q)_{\infty}$ d'\^etre finie. Ainsi, lorsque $z$ parcourt $\mathbb{C}_{q,m}$ on obtiendra  un ensemble d'\'{e}tats coh\'{e}rents $\left|z,m,q\right\rangle$ tels que l'int\'{e}grat-ion par rapport \`a une mesure convenable  des op\'{e}rateurs de rang un $\left|z,m,q\right\rangle\left\langle z,m,q \right|$  donne lieu \`a l'\'{e}galit\'{e}
\begin{equation}\label{normalisation}
\int_{\mathbb{C}_{q,m}} |z,m,q\rangle\langle z,m,q|d\nu_{m,q}(z)=1_{\mathcal{H}_q}
\end{equation}
laquelle traduit la r\'esolution de l'identit\'e $1_{\mathcal{H}_q}$ de l'espace de Hilbert $\mathcal{H}_q:=L^2(\mathcal{I}_q,d\xi)$ suppos\'e abriter  les \'etats quantiques d'un oscillateur harmonique $q$-d\'eform\'e. Pr\'{e}cisement, la mesure dans  \eqref{normalisation} s'\'ecrit 
\begin{equation}
d\nu_{m,q}(z)=\mathcal{N}_{m,q}(z) d\mu_q(z)
\end{equation}
avec
\begin{equation}\label{dmuq}
d\mu_q(z)=\sum_{j\geq 0}\frac{q^j(q;q)_{\infty}}{(q;q)_j}d\mu_j(z),
\end{equation}
$d\mu_j(z)$ \'etant la mesure de Lebesgue sur le cercle de rayon $\rho_j=q^{\tfrac{1}{2}j}(1-q)^{-1/2}$. Sachant que la superposition dans \eqref{CS1} servira,  \`a une racine carr\'ee de $\mathcal{N}_{m,q}(z)$ multiplicative pr\`es, de noyau pour la transformation qu'on se propose de construire, on aura besoin d'une forme compacte pour la somme dans  \eqref{CS1}. D'o\`u l'int\'er\^et du r\'esultat ci-dessous.\\ \\
\textbf{Proposition 1.1.} \label{propo}\textit{Soient $m\in \mathbb{Z}_+$ et $0<q<1$. Alors pour chaque $z\in\mathbb{C}_{q,m}$ fix\'e, la fonction d'onde de l'\'etat coh\'erent \eqref{CS1} a pour expression}
\begin{eqnarray}\label{SalCHi}
\langle \xi|z,m,q \rangle&=&(-1)^m\left(\frac{\sqrt{2}\;\omega_q(\sqrt{2}\xi)}{q^m(q;q)_m\mathcal{N}_{m,q}(z\bar{z})}\right)^{\tfrac{1}{2}} \frac{1}{\Big|(ze^{\,i\arccos(\xi\sqrt{\tfrac{1-q}{2}})}\sqrt{\frac{1-q}{q^m}};q)_{\infty}\Big|^2}\cr
&\times& Q_m\left(\sqrt{\tfrac{1-q}{2}}\xi;\sqrt{\tfrac{1-q}{q^m}}z,\sqrt{\tfrac{1-q}{q^m}}q\bar{z}|q\right)
\end{eqnarray}\medskip\\
\textit{en termes des polyn\^omes d'AL-Salam-Chihara} $Q_m$  \textit{et ce pour tout $\xi\in \mathbb{R}.$}\medskip\\
Rappelons que les polyn\^omes qui figurent dans \eqref{SalCHi} se d\'efinissent \`a l'aide de la s\'erie $q$-hyperg\'eom\'etrique ${}_3\phi_2$ comme suit   (\cite{KS}, p.80) :
\begin{equation}
Q_m(x;a,b|q):=\frac{(ab;q)_m}{a^m}\setlength\arraycolsep{1pt}
{}_3 \phi_2\left(\begin{matrix}q^{-m},ae^{i\theta},ae^{-i\theta} \\ ab,0 \end{matrix}\left|q;q\right.\right)
\end{equation}
o\`u $x=$cos$\, \theta$.\medskip

Dans le cas particulier o\`u $m=0$, on obtient par un remplacement direct ce qui suit. \medskip\\ 
\textbf{Corollaire 1.1} \textit{Soit $0<q<1$. La fonction d'onde \eqref{CS1} qui correspond au niveau } $m=0$ \textit{est de la forme }
\begin{equation}\label{CS2}
\langle \xi|z,0,q\rangle = \left(\frac{\sqrt{2}\;\omega_q(\sqrt{2}\xi)}{e_q(z\bar{z})}\right)^{\tfrac{1}{2}}\displaystyle\prod_{k\geq 0}\frac{1}{\left(1-\sqrt{2}\bar{z}\xi q^k(1-q)+\bar{z}^2q^{2k}(1-q)\right)},
\end{equation}
\textit{o\`u $z\in \mathbb{C}_{q,m}$ \'etant fix\'e et $\xi\in \mathbb{R}$}.\medskip\\
Un tel \'etat coh\'erent ne pourra \^etre normalis\'e que si la variable d'indexation $z$ reste dans le domaine  $\mathbb{C}_q:=\mathbb{C}_{q,0}=\{z\in \mathbb{C},(1-q)z\bar{z}<1\}$ qui n'est autre que celui de la convergence de la s\'erie $e_q(z\bar{z})$. Ici,
\begin{equation}
e_q(u)=\displaystyle\sum_{k\geq 0}\frac{u^k}{[k]_q!},
\end{equation}
o\`u $(q;q)_k=[k]_q!(1-q)^k.$\medskip

 En restant toujours dans le cas $m=0$, notons par $\mathcal{A}^2(\mathbb{C}_q)$ le compl\'et\'e de l'espace des fonctions holomorphes sur $\mathbb{C}_q$, muni du produit  scalaire
\begin{equation*}
\langle\varphi,\phi\rangle=\int_{\mathbb{C}_q}{\varphi(z)}\overline{\phi(z)}d\mu_q(z)
\end{equation*}
o\`u $d\mu_q$ d\'esigne la mesure donn\'ee ci-dessus par \eqref{dmuq}, il est bien connu que cette mesure est unique \`a rendre le syst\`eme  $h_j^{0,q}(z)=([j]_q!)^{-1/2}z^j$  une base orthonormale de l'espace $\mathcal{A}^2(\mathbb{C}_q)$ qui n'est autre que l'espace de Arik-Coon \cite{ACO}. Dans cet espace, l'op\'erateur d'annulation des bosons se repr\'esente par un op\'erateur aux $q$-diff\'erence qui tendra vers l'op\'erateur de d\'erivation quand $q\rightarrow 1$. A la m\^eme limite la mesure $d\mu_q$ deviendra la mesure Gaussienne sur $\mathbb{C}$.\medskip

Dans le cas o\`u  $m$ est non nul,  par un calcul direct, nous \'etablissons la relation de chevauchement entre deux \'etats coh\'erents quelconques.\medskip\\
\textbf{Proposition 1.2.} \textit{Soient $m\in\mathbb{Z}_+$, $q\in]0,1[$ et $z,w\in \mathbb{C}_{q,m}$. Alors, on a}
\begin{equation}\label{kern1}
\langle z,m,q| w,m,q\rangle = \frac{q^{-m}(q^{-m+1}(1-q)z\bar{z};q)_{\infty}}{ ( q^{-m}(1-q)w\bar{z},q(1-q)z\bar{z};q)_\infty}\frac{(\frac{\bar{w}}{\bar{z}} \,q,q)_m}{(q;q)_m}
{}_3\phi_2\left(\begin{array}{c}q^{-m},  q^{-m}(1-q)w\bar{z}, \frac{\bar{w}}{\bar{z}}\\
 q^{-m+1}(1-q)z\bar{z},\frac{\bar{z}}{\bar{w}}q^{-m}\end{array}\Big|q;q \right).
\end{equation}
\medskip\\ D'une autre part, comme il a \'et\'e mentionn\'e ci dessous  cette derni\`ere \'egalit\'e permet de retrouver par un calcul direct le facteur de normalisation  \eqref{normfac} en posant $z=w$. D'autre part, si l'on introduit la fonction 
\begin{equation}\label{repkern}
K_{m,q}(z,w):= \langle z,m,q| w,m,q\rangle
\end{equation}
et que l'on d\'efinisse l'espace $\mathcal{A}^2(\mathbb{C}_{q,m})$ des fonctions de carr\'e int\'egrables sur $\mathbb{C}_{q,m}$ par rapport \`a la mesure $d\mu_q$ dans  \eqref{dmuq} et ayant pour noyau reproduisant la fonction \eqref{repkern} alors on peut \'enoncer ce qui suit.     
\medskip \\
\textbf{Th\'eor\`eme 1.} \textit{Pour $m\in\mathbb{Z}_+$ et $q\in]0,1[$.  La transformation issue des \'etats coh\'erents \eqref{CS1} est l'isom\'etrie $\mathcal{B}_{m}^q:\mathcal{H}_q\longrightarrow \mathcal{A}^2(\mathbb{C}_{q,m})$, d\'efinie  par }
\begin{equation}
\mathcal{B}_{m}^q[f](z)=\frac{(-1)^m}{\sqrt{q^m(q;q)_m}}\displaystyle\int_{\mathbb{R}}\frac{Q_m\left(\sqrt{\tfrac{1-q}{2}}\xi;\sqrt{\tfrac{1-q}{q^m}}z,\sqrt{\tfrac{1-q}{q^m}}q\bar{z}|q\right)}{\Big|(ze^{i\,\arccos(\xi\sqrt{\tfrac{1-q}{2}})}\sqrt{\frac{1-q}{q^m}};q)_{\infty}\Big|^2} 
\sqrt{\sqrt{2}\;\omega_q(\sqrt{2}\xi)}f(\xi)d\xi
\end{equation}
\textit{en tout point $z\in \mathbb{C}_{q,m}$.}
\medskip\\
\textbf{D\'efinition 2.1.} \textit{L'isom\'etrie $\mathcal{B}_m^q$ est appel\'ee la $q$-d\'eform\'e de la transformation de Bargmann polyanalytique.}\medskip\\
\textbf{Corollaire 1.2.} \textit{La  $q$-d\'eform\'ee de la transformation de Bargmann qui correspond au cas analytique} $m=0$ \textit{est l'isom\'etrie $\mathcal{B}_0^q:\mathcal{H}_q\longrightarrow\mathcal{A}^2(\mathbb{C}_q)$, d\'efinie  par }
\begin{equation}\label{B0q}
\mathcal{B}_q^0[f](z)=\displaystyle\int_{\mathbb{R}}\left(\displaystyle\prod_{k\geq 0}\frac{1}{\left(1-\sqrt{2}{z}\xi q^k(1-q)+{z}^2q^{2k}(1-q)\right)}\right)\sqrt{\sqrt{2}\;\omega_q(\sqrt{2}\xi)}f(\xi)d\xi.
\end{equation}
Bien entendu, comme il fallait s'y attendre l'\'equation \eqref{B0q} permet de retrouver la transformation de Bargmann classique \eqref{barcl}  lorsqu'on fait tendre $q$ vers 1. Cela se justifie par les limites 
\begin{equation}
\displaystyle\lim_{q\rightarrow 1}\displaystyle\prod_{k\geq 0}\frac{1}{\left(1-\sqrt{2}{z}\xi q^k(1-q)+{z}^2q^{2k}(1-q)\right)}=e^{\sqrt{2}\xi z-\tfrac{1}{2}z^2}
\end{equation}
et
\begin{equation}
\displaystyle\lim_{q\rightarrow 1} \omega_q(u)=\frac{1}{\sqrt{2\pi}}e^{-\tfrac{1}{2}u^2}
\end{equation}
qui furent \'etablies dans (\cite{MDG}, p. 381-382) avec la notation $\omega=\tilde{v}$.\medskip\\
\textbf{Remarque 1.1.} Dans le cas $m=0$, si l'on pose
\begin{equation}
z=2\alpha\:,\:\xi=\frac{\sqrt{2}}{\sqrt{1-q}}\cos\theta
\end{equation}
et que l'on  d\'esigne par $\phi_0(\theta)$ la fonction telle que 
\begin{equation}
\left(\phi_0(\theta)\right)^2=\sqrt{2}\:\omega_q\left(\frac{\sqrt{2}}{\sqrt{1-q}}\cos\theta\right)
\end{equation}
on peut alors s'assurer que l'expression ainsi obtenue \`a partir de \eqref{CS2} 
\begin{equation}
({e_q(4\alpha\bar{\alpha})})^{-\tfrac{1}{2}}\left[\;\:\phi_0(\theta)\:\frac{1}{(2\alpha e^{i\theta};q)_{\infty}(2\alpha e^{-i\theta};q)_{\infty}}\right]
\end{equation}  
ne diff\`ere que par le facteur de normalisation $\left(\mathcal{N}_{0,q}(2\alpha) \right)^{-1/2}$ de celle des \'etats coh\'erents que Odake et Sasaki ont construit \textquotedblleft \textit{\`a la Glauber} \textquotedblright pour l'oscillateur harmonique $q$-d\'eform\'e (\cite{OS}, p.144, Eq.52).\medskip
\section{Esquisse de d\'emonstration des r\'esultats}

Pour \'etablir l'\'enonc\'e de la Proposition {\bf 1.1} dont on s'est servi pour conclure le r\'esultat principal, on commence par remplacer le quantit\'es $\varphi_j^q(\xi)$ et $h_j^{m,q}(z)$ qui interviennent dans \eqref{CS1}  par leurs expressions dans  \eqref{Hmj} et \eqref{fjq123} respectivement. Cela conduit \`a la sommation suivante 
\begin{eqnarray}\label{So31}
|z,m,q\rangle&=&(\mathcal{N}_{m,q}(z\bar{z}))^{-\tfrac{1}{2}}\sum_{j\geq 0} \frac{(-1)^{m\wedge j} (q;q)_{m \vee j}q^{\binom{m\wedge j}{2}}\sqrt{1-q}^{|m-j|}|z|^{|m-j|}e^{-i(m-j)\theta}}{(q;q)_{|m-j|}\sqrt{q^{mj}(q;q)_m(q;q)_j}}\cr &\times& P_{m\wedge j}\left( (1-q)z\bar{z};q^{|m-j|}|q\right)
  \sqrt{\frac{\sqrt{2}\,\omega_q(\sqrt{2}\xi)}{(q;q)_j}}H_j\left(\sqrt{\frac{1-q}{2}}\xi|q\right).
\end{eqnarray}
Dor\'enavant, on s'int\'eressera qu'\`a la somme qui appara\^it dans \eqref{So31} et qu'on re\'ecrit sous la forme de deux morceaux  $S_{(<\infty)}\left(z,q,m;\xi\right)+S_{(\infty)}\left(z,q,m;\xi\right)$ o\`u 
 \begin{eqnarray}
 \mathcal{S}_{(<\infty)}^q(m,z;\xi)&=& \displaystyle\sum_{j=0}^{m-1} \frac{(-1)^{ j} (q;q)_{m }q^{\binom{ j}{2}}\sqrt{1-q}^{m-j}\bar{z}^{m-j}}{(q;q)_{m-j}\sqrt{q^{mj}(q;q)_m(q;q)_j}}P_{j}\left( (1-q)z\bar{z};q^{m-j}|q\right) \varphi_j^q(\xi)\cr
 &-&  \displaystyle\sum_{j=0}^{m-1} \frac{(-1)^{m} (q;q)_{ j}q^{\binom{m}{2}}\sqrt{1-q}^{j-m}z^{j-m}}{(q;q)_{j-m}\sqrt{q^{mj}(q;q)_m(q;q)_j}}P_{m}\left( (1-q)z\bar{z};q^{j-m}|q\right) \varphi_j^q(\xi),
 \end{eqnarray}
et
\begin{eqnarray}
S_{(\infty)}\left(z,q,m;\xi\right)&=&\frac{(-1)^mq^{({}^m_2)}(z\sqrt{1-q})^{-m}\sqrt{\sqrt{2}\omega_q(\sqrt{2}\xi)}}{\sqrt{(q;q)_m}}\sum_{j=0}^\infty\frac{\sqrt{(q;q)_j}\;Y^j}{(q;q)_{j-m}}\cr
&\times&P_m((1-q)z\bar{z}; q^{j-m}|q)\frac{ H_j(X|q)}{\sqrt{(q;q)_j}}
\end{eqnarray}
avec $Y=z\sqrt{\frac{1-q}{q^m}}$ et $X=\sqrt{\frac{1-q}{2}} \xi$ . En faisant appel \`a l'identit\'e (\cite{MoCa},  p.3) :
\begin{equation}
 P_n(x;q^{-N}|q)=x^N(-1)^{-N}q^{\frac{N(N+1-2n)}{2}}\frac{(q^{N+1};q)_{n-N}}{(q^{1-N};q)_n} P_{n-N}(x;q^N|q)
 \end{equation}
 que v\'erifient les p\^olynomes de Wall et ce pour les param\`etres $N=j-m,n=j$ et $x=(1-q)z\bar{z}$, on peut se convaincre que la somme finie vaut effectivement z\'ero. Concernant la somme infinie, on re\'ecrit le polyn\^ome de Wall comme suit   
 \begin{equation}
 P_m(x;q^{j-m}|q)=\frac{1}{(q^{j-m+1};q)_m} \sum_{r,k=0}^{\infty}(q^{-m};q)_{r+k}\frac{(qx)^r}{(q;q)_r}\frac{(q^{j+1})^k}{(q;q)_k}.
 \end{equation}
Cette formule peut \^etre d\'eduite de la formule g\'en\'eratrice
\begin{equation}\label{So32}
\sum_{r,k=0}^{+\infty} \mathcal{P}_{r+k}(a,b)\frac{t^r}{(q;q)_r}\frac{s^k}{(q;q)_k}=\frac{(bs;q)_{\infty}}{(as;q)_{\infty}}\setlength\arraycolsep{1pt}
{}_2 \phi_1\left(\begin{matrix}b/a,0 \\ bs \end{matrix}\left|q;at\right.\right)
\end{equation} 
des polyn\^omes de Cauchy $\mathcal{P}_n$ et de \eqref{Wall1} qui exprime les polyn\^omes  $P_n$ de Wall en termes de $_2\phi_1$. Cela donne 
\begin{eqnarray}
S_{(\infty)}\left(z,q,m;\xi\right)&=&\frac{(-1)^mq^{({}^m_2)}(z\sqrt{1-q})^{-m}\sqrt{\sqrt{2}\omega_q(\sqrt{2}\xi)}}{\sqrt{(q;q)_m}} \sum_{j,r,k=0}^{\infty}\frac{Y^j}{(q;q)_{j-m}(q^{j-m+1};q)_m}\cr
&\times&(q^{-m};q)_{r+k}\frac{(q\xi)^r}{(q;q)_r}\frac{(q^{j+1})^k}{(q;q)_k}
H_j(X|q)
\end{eqnarray}
que l'on ordonne 
\begin{equation*}
S_{(\infty)}\left(z,q,m;\xi\right)=\frac{(-1)^mq^{({}^m_2)}(z\sqrt{1-q})^{-m}\sqrt{\sqrt{2}\omega_q(\sqrt{2}\xi)}}{\sqrt{(q;q)_m}} \sum_{r,k=0}^{\infty} (q^{-m};q)_{r+k} \frac{(q\xi)^r}{(q;q)_r}\frac{q^k}{(q;q)_k} 
\end{equation*}
\begin{equation}\label{So33}
 \times\sum_{j=0}^{\infty} \frac{(Yq^k)^j}{(q;q)_j}H_j(X|q)
\end{equation}
de sorte \`a pouvoir utiliser la fonction g\'en\'eratrice des polyn\^omes $q$-Hermite continus. Ainsi, l'\'equation \eqref{So33} devient 
\begin{eqnarray}
S_{(\infty)}\left(z,q,m;\xi\right)&=&\frac{(-1)^mq^{({}^m_2)}(z\sqrt{1-q})^{-m}\sqrt{\sqrt{2}\omega_q(\sqrt{2}\xi)}}{\sqrt{(q;q)_m}} \frac{1}{|(e^{i\theta}Yq;q)_{\infty}|^2} \cr &\times& \sum_{k=0}^{\infty}\frac{(e^{i\theta}Y;q)_k(e^{-i\theta}Y;q)_k}{(q;q)_k}q^k\sum_{r=0}^{\infty} \frac{(q\xi)^r}{(q;q)_r}(q^{-m};q)_{r+k} .
\end{eqnarray}
A ce stade on applique l'identit\'e (\cite{SA}, p.2):

\begin{equation}
\sum_{n\geq 0} \frac{(\lambda;q)_{m+n}}{(q;q)_n}t^n=\frac{(\lambda;q)_m}{(\lambda t;q)_m}\;\frac{(\lambda t;q)_{\infty}}{(t;q)_{\infty}}\:,\qquad |t|<1,\qquad |q|<1,
\end{equation}
dans le cas des param\`etres $\lambda=q^{-m}$ et $t=q\xi$. Ce qui nous fait aboutir \`a la forme 
\begin{eqnarray}
S_{(\infty)}\left(z,q,m;\xi\right) &=&\frac{(-1)^mq^{({}^m_2)}(z\sqrt{1-q})^{-m}\sqrt{\sqrt{2}\omega_q(\sqrt{2}\xi)}}{\sqrt{(q;q)_m}} \frac{(q^{1-m}\xi;q)_{\infty}}{(\xi q;q)_{\infty}|(e^{i\theta}Y;q)_{\infty}|^2} \cr
&\times& \sum_{k=0}^{\infty}\frac{(q^{-m};q)_k(e^{i\theta}Y;q)_k(e^{-i\theta}Y;q)_k}{(q^{1-m}\xi;q)_k}\,\frac{q^k}{(q;q)_k}.
\end{eqnarray}
On reconnait la somme qui figure dans la derni\`ere expression comme \'etant la s\'erie $q$-hyperg\'eom\'etrique

\begin{equation}
{}_3 \phi_2\left(\begin{matrix}q^{-m},Ye^{i\theta},Ye^{-i\theta} \\ q\sqrt{\frac{1-q}{q^m}}\bar{z},0 \end{matrix}\left|q;q\right.\right)\:,
\end{equation}
qui au fait contient un nombre fini de termes et d\'efinit le polyn\^ome d'Al-Salam-Chihara
\begin{equation}
\frac{Y^m}{(q^{1-m}\xi;q)_m} Q_m(X;Y,q\sqrt{\frac{1-q}{q^m}}\bar{z}|q).
\end{equation}
Enfin, en tenant compte de tous les pr\'efacteurs qui ont apparu lors du calcul dans chacune des \'etapes pr\'ec\'edentes, on aboutit \`a la forme de la fonction d'onde de l'\'etat coh\'erent d\'efini dans \eqref{SalCHi} moyennant le facteur de normalisation $\mathcal{N}_{m,q}(z)$. Ce dernier s'obtient on posant $z=w$ dans la relation \eqref{kern1}. Finalement on dispose d'une expression compact de la fonction d'onde permettant de conclure le r\'esultat du th\'eor\`eme.

\end{document}